\documentclass[12pt,reqno]{amsart}
\usepackage[headings]{fullpage}
\usepackage{amssymb,amsmath,amscd,bbm,tikz}
\usepackage{graphicx}
\usepackage{texdraw}
\usepackage{pb-diagram}
\usepackage[all,cmtip]{xy}
\usepackage{url}
\usepackage[bookmarks=true,%
    colorlinks=true,%
    linkcolor=blue,%
    citecolor=blue,%
    filecolor=blue,%
    menucolor=blue,%
    urlcolor=blue,%
    breaklinks=true]{hyperref}
\usepackage{slashed}    
\usepackage{verbatim}
\usepackage[normalem]{ulem}    
\usepackage{cancel}            
\usepackage{lmodern}           

\newtheorem{theorem}{Theorem}[section]
\theoremstyle{definition}

\newtheorem{corollary}[theorem]{Corollary}
\newtheorem{conjecture}[theorem]{Conjecture}

\def\BW{\mathbbm W}

\def\BZ{\mathbbm Z}
\def\BQ{\mathbbm Q}
\def\BR{\mathbbm R}
\def\BC{\mathbbm C}

\def\BW{\mathbbm W}

\def\calE{\mathcal E}
\def\calT{\mathcal T}

\def\calM{\mathcal M}

\def\calO{\mathcal O}
\def\s{\sigma}

\def\SL{\mathrm{SL}}

\def\pt{\partial}
\def\ID{I_{\Delta}}
\def\tq{\tilde{q}}

\def\a{\alpha}
\def\b{\beta}

\def\PSL{\mathrm{PSL}}

\def\be{\begin{equation}}
\def\ee{\end{equation}}
\def\ID{I_{\Delta}}

\usepackage{accents}

\def\hz{\widehat{z}}
\def\Irot{I^{\mathrm{rot}}}

\def\Ahat{\widehat{A}}

\def\wh{\widehat}

\newcommand{\mb}{\mathbf}

\makeatletter

\renewcommand\thepart{\@Roman\c@part}%
\renewcommand\part{%
   \if@noskipsec \leavevmode \fi
   \par
   \addvspace{6.7ex}%
   \@afterindentfalse
   \secdef\@part\@spart}
\def\@part[#1]#2{%
    \ifnum \c@secnumdepth >\m@ne
      \refstepcounter{part}%
      \addcontentsline{toc}{part}{Part~\thepart.\ #1}%
    \else
      \addcontentsline{toc}{part}{#1}%
    \fi
    {\parindent \z@ \raggedright
     \interlinepenalty \@M
     \normalfont
     \ifnum \c@secnumdepth >\m@ne
       \centering\large\scshape \partname~\thepart.%
       \hspace{1ex}%
     \fi%
     \large\scshape #2%
     \markboth{}{}\par}%
    \nobreak
    \vskip 4.7ex
    \@afterheading}
  \def\@spart#1{
  \refstepcounter{part}%
  \addcontentsline{toc}{part}{#1}%
    {\parindent \z@ \raggedright
     \interlinepenalty \@M
     \normalfont
     \centering\large\scshape #1\par}%
     \nobreak
     \vskip 4.7ex
     \@afterheading}
\renewcommand*\l@part[2]{%
  \ifnum \c@tocdepth >-2\relax
    \addpenalty\@secpenalty
    \addvspace{0.75em \@plus\p@}%
    \begingroup
      \parindent \z@ \rightskip \@pnumwidth
      \parfillskip -\@pnumwidth
      {\leavevmode
       \normalsize \bfseries #1\hfil \hb@xt@\@pnumwidth{\hss #2}}\par
       \nobreak
       \if@compatibility
         \global\@nobreaktrue
         \everypar{\global\@nobreakfalse\everypar{}}%
      \fi
    \endgroup
  \fi}

\def\l@subsection{\@tocline{2}{0pt}{2pc}{6pc}{}}
\makeatother

\begin{document}
\title[The descendants of the 3d-index]{
       The descendants of the 3d-index}
\author{Zhihao Duan}
\address{Korea Institute for Advanced Study,
85 Hoegiro, Dongdaemun-gu, Seoul 02455, Korea}
\email{xduanz@gmail.com}
\author{Stavros Garoufalidis}
\address{
  International Center for Mathematics, Department of Mathematics \\
  Southern University of Science and Technology \\
  Shenzhen, China \newline
  {\tt \url{http://people.mpim-bonn.mpg.de/stavros}}}
\email{stavros@mpim-bonn.mpg.de}
\author{Jie Gu}
\address{School of Physics and Shing-Tung Yau Center\\
Southeast University, Nanjing 210096, China}
\email{eij.ug.phys@gmail.com}

\begin{flushright}
{\tt\normalsize KIAS-P22086}\\
\end{flushright}


\date{1 January 2023}

\begin{abstract}
  In the study of 3d-3d correspondence occurs a natural $q$-Weyl algebra associated
  to an ideal triangulation of a 3-manifold with torus boundary components, and a
  module of it. We study the action of this module on the (rotated) 3d-index of
  Dimofte--Gaiotto--Gukov and we conjecture some structural properties:
  bilinear factorization in terms of holomorphic blocks, pair of linear
  $q$-difference equations, the determination of the 3d-index in terms of a
  finite size matrix of rational functions and the asymptotic expansion of the
  $q$-series as $q$ tends to 1 to all orders. We illustrate our conjectures with
  computations for the case of the three simplest hyperbolic knots. 
\end{abstract}

\maketitle

{\footnotesize
\tableofcontents
}


\section{Introduction}
\label{sec.intro}

\subsection{The 3D-index and the state-integral}
  
Topological invariants of ideally triangulated 3-manifolds appeared in
mathematical physics in relation to complex Chern--Simons
theory~\cite{AK:TQFT} and its extension in the 3d-3d
correspondence~\cite{Dimofte:riemann,Gang:aspects}. Two of the
best-known such invariants are the state-integrals of
Andersen--Kashaev~\cite{AK:TQFT}, which are analytic functions on
$\BC\setminus (-\infty,0]$, and the 3D-index of
Dimofte--Gaiotto--Gukov~\cite{DGG1,DGG2}, which is a collection of
$q$-series with integer coefficients parametrized by the integer
homology of the boundary of a 3-manifold. Although the state-integrals
and the 3D-index are different looking functions, they are closely
related on the mathematics side through the theory of holomorphic
quantum modular forms developed by Zagier and the second
author~\cite{GZ:kashaev,GZ:qseries}, and on the physics side through
the above mentioned 3d-3d correspondence.

The state-integrals and the 3D-index share many common features, stemming
from the fact that on the physics side, under the 3d-3d correspondence
\cite{DimofteGukovHollands2011,TerashimaYamazaki2011,DGG2,
  DimofteGabellaGoncharov2016} (see \cite{Dimofte:3dsuper} for a review)
become the invariants of the
dual 3d $N=2$ superconformal field theory on respectively
$S^3$ and $S^1 \times S^2$, both of which can be obtained by gluing
two copies of $D^2\times S^1$ together.

On the mathematics side, both invariants are defined using
combinatorial data of ideal triangulations of 3-manifolds whose local
weights (namely the Faddeev quantum dilogarithm function, and the
tetrahedron index, respectively) satisfy the same linear
$q$-difference equations, whereas the invariants themselves are given
by an integration/summation over variables associated to each
tetrahedron.

A common feature to both invariants is their conjectured bilinear
factorization in terms of the same holomorphic blocks $H(q)$, the
latter being $q$-hypergeometric series defined for $|q| \neq 1$. This
leads to bilinear expressions for the state-integral in terms of
$H(q)$ times $H(\tq)$ (where $q=e^{2 \pi i \tau}$ and
$\tq=e^{-2 \pi i/\tau}$) and bilinear expressions for the 3D-index in
terms of $H(q)$ times $H(q^{-1})$. This factorization is well-known
in the physics literature~\cite{Beem} and interpreted as
partition function of the dual 3d superconformal field theory on
$D^2\times S^1$. They are also partially known for some examples of
3-manifolds in~\cite{GK:qseries, GZ:qseries}.  We emphasize, however,
that the bilinear factorization of state-integrals and of the 3D-index
is conjectural, and so is the existence of the suitably normalized
holomorphic blocks.

Another common feature to state-integrals and the 3D-index is that
they are given by integrals/lattice sums where the integrand/summand
has a common annihilating ideal. This implies that both
state-integrals and the rotated 3D-index satisfy a pair of linear
$q$-difference equations which are in fact conjectured to be
identical, and equal to the homogeneous part of the linear
$q$-difference equation for the colored Jones polynomial of a
knot~\cite{GL:qholo}. The conjectured common linear $q$-difference
equations for state-integrals and for the 3D-index would also be a
consequence of their common holomorphic block
factorization. In physics these linear $q$-difference
equations are interpreted as Ward identities of Wilson-'t Hooft line
operators in the dual 3d superconformal field theory \cite{DGG1,DGG2}.

\subsection{Descendants}
\label{sub.descendants}

Descendants appeared recently as computable, exponentially small corrections
to the asymptotics of the Kashaev invariant of a knot, refining the Volume
Conjecture to all orders in perturbation theory to a Quantum Modularity
Conjecture~\cite{GZ:kashaev}. One of the discoveries was that the Kashaev invariant of
a knot is a distinguished $(\s_0,\s_1)$-entry in a square matrix of knot invariants
at roots of unity. The rows and columns of the matrix are parametrized with
boundary-parabolic $\PSL_2(\BC)$-representations, with $\s_0$ denoting the trivial
representation and $\s_1$ denoting the geometric representation of a hyperbolic
knot complement. The above mentioned matrix has remarkable algebraic, analytic and
arithmetic properties explained in detail in Section 5 of~\cite{GZ:kashaev},
and given explicitly for the $4_1$ and $5_2$ knots in Sections 7.1 and 7.2 of i.b.i.d.
The rows of the matrix are supposed to be $\BQ(q^{1/2})$-linear combinations of
fundamental solutions to a linear $q$-difference equation (homogeneous for all but
the first row), thus the elements in each row are supposed to be descendants of
each other. 
Although the existence of such a matrix is conjectured, its top row was defined
in~\cite{GK:desc} for all knots in terms of the descendant Kashaev invariants of a knot. 

The above mentioned matrix has three known realizations, one as functions at roots
of unity mentioned above, a second as a matrix of Borel summable asymptotic series
and a third as a matrix of $q^{1/2}$-series. The idea of descendants can be extended
to the matrix of asymptotic series (whose first column are simply the vector of
asymptotic series of the perturbative Chern--Simons theory at a $\PSL_2(\BC)$-flat
connection, and the remaining columns being descendants of the first column)
as well as to a matrix of $q$-series. This extension was done for the case of the
$4_1$ and $5_2$ knots by Mari\~{n}o and two of the authors~\cite[Eqn.(13),App.A]{GGM},
with the later addition of the trivial $\PSL_2(\BC)$-representation
in~\cite[Sec.2.2,Sec.4.1]{GGMW:trivial}.

To summarize, descendants are supposed to be the $\BQ(q^{1/2})$-span of a fundamental
solution to a linear $q$-difference equation associated to the quantum invariants.
It is becoming clear that this span is a fundamental quantum invariant of 3-manifolds,
and we want to present further evidence for this using as an example an important
quantum invariant, namely the 3D-index. 

\subsection{Our conjectures}
\label{sub.summary}

A detailed study of the 3D-index of a 3-manifold with torus boundary and its
structural properties, namely holomorphic block factorization, linear $q$-difference
equations, computations and asymptotics was recently done in~\cite{GW:periods}.

The goal of the present paper is to extend the properties of the
3D-index by allowing observables, line operators, defects,
descendants, all being synonymous names for the same object. On the
topological side, an observable is a knot $L$ in a 3-manifold $M$,
where in the case of interest to us, $M=S^3\setminus K$ is the
complement of a knot in $S^3$. On the algebra side, the conjectural
3d-quantum trace map sends a knot $L \subset S^3\setminus K$ to an
element $\calO$ of a module over a $q$-Weyl algebra associated to an
ideal triangulation $\calT$ of $M$. We will postpone the description
of the 3d-quantum trace map to a subsequent publication. Now $\calO$
acts on the integrand/summand of the state-integral/3D-index, and by
integrating/summing one obtains a state-integral/3D-index with
insertion $\calO$. On the physics side, $\calO$ becomes a
line-operator supported on a line $\gamma$ in the dual 3d $N=2$
superconformal field theory $T_2[M]$ under the 3d-3d
correspondence~\cite{DimofteGukovHollands2011,TerashimaYamazaki2011,
  DGG2,DimofteGabellaGoncharov2016}.
The 3d-3d correspondence can be understood as a consequence of
compactifying 6d $N=2$ $A_1$ superconformal field theory on the three
manifold $M$ and on $\BR^3$ with topological twist along
$M$.  The 6d theory has surface operators which can be supported on
$L\times \gamma$, giving rise to the correspondence between the defect
$L$ in $M$ and the line-operator on $\gamma \subset \BR^3$ in $T_2[M]$
\cite{DGG1,DGG2}. Our goal
is to study the structural properties of the rotated, inserted,
3D-index $\Irot_{\calT,\calO}(q)$. Although this is a $\BZ \times \BZ$
matrix, we will see that it is determined from the uninserted rotated
3D-index $\Irot_\calT(q)$ in terms of a pair of linear $q$-difference
equations and a finite size invertible matrix with coefficients in the
field $\BQ(q^{1/2})$; see Conjectures~\ref{conj.1} and~\ref{conj.2}
below, illustrated by examples in Section~\ref{sec.examples}.

We emphasize that our paper concerns conjectural structural properties of
topological invariants, such as the rotated inserted 3D-index, and not
mathematical proofs. Nevertheless the structure of these invariants is rich,
and leads to startling predictions and numerical conformations
(see eg. Equation~\eqref{cancel} below). 


\section{Algebras of 3-dimensional ideal triangulations}
\label{sec.alberas}

We recall here a $q$-Weyl algebra associated to an
ideal triangulation $\calT$ which was first considered by Dimofte on the context of
the 3d-3d correspondence, and it was introduced as an attempt to quantize the
$\SL_2(\BC)$-character variety of an ideally triangulated 3-manifold $M$ using
the symplectic structure of the Neumann--Zagier matrices, and following the ideas
of Hamiltonian reduction of symplectic phase-spaces~\cite{Dimofte:riemann}.
Similar ideas appeared in subsequent work~\cite{Gang:aspects}. 

We fix an ideal triangulation
$\calT$ of $M$ with $N$ ideal tetrahedra. This defines a $q$-Weyl algebra
$\BW_q(\calT)=\BQ(q)\langle \hz_j, \hz_j' \, | \, j=1,\dots,N \rangle$ of Laurent
variables $\hz_j, \hz_j'$ that commute except in the following instance
$\hz_j \hz_j'=q \hz_j' \hz_j$ for $j=1,\dots,N$. A more symmetric way is to
introduce three invertible variables $\hz, \hz',\hz^{''}$ which satisfy the relations
\be
\label{zrels}
\hz \hz' = q \hz' \hz, \qquad \hz' \hz'' = q \hz'' \hz', \qquad
\hz'' \hz = q \hz \hz'', \qquad \hz \hz' \hz'' = -q 
\ee
(hence  $\hz \hz' \hz''$ is in the center and it is invariant under cyclic
permutations), and then $\BW_q(\calT)$ is simply the tensor product of one algebra
per tetrahedron.
The combinatorics of the edge-gluing equations of $M$ have
symplectic properties discovered by Neumann--Zagier~\cite{NZ,N:comb}. Using
those properties, Dimofte~\cite{Dimofte:riemann} and later
Gang et al~\cite{Gang:aspects} (see also~\cite[Eqn.(10)]{Agarwal})
consider the quotient \be
\label{strangequotient}
\calM(\calT) = \BW_q(\calT)/(\BW_q(\calT)(\text{Lagrangians})
+(\text{edge equations})\BW_q(\calT))
\ee
of $\BW_q(\calT)$ by the left $\BW_q(\calT)$-ideal generated by the Lagrangian
equations
\be
\label{3lag}
\hz'^{-1} + \hz -1=0, \qquad (\hz^{''})^{-1} + \hz' -1=0, \qquad
\hz^{-1} + \hz^{''} -1=0 
\ee
(one per each tetrahedron) plus the right ideal generated by the
edge equations. This strange quotient $\calM(\calT)$, which is no longer a
module over a $q$-Weyl algebra, but only a $\BQ(q^{1/2})$-vector space is a natural
object that indeed annihilates the rotated 3D-index as we will see shortly.


\section{The rotated 3D-index and its descendants}
\label{sec.review}

\subsection{Definition}
\label{sub.def}

For simplicity, in the paper we will focus on the action of the quantum torus
$\BW_q(\calT)$ on the 3D-index $I_\calT$, and in fact in its rotated form $\Irot_\calT$
explained to us by Tudor Dimofte and studied extensively in~\cite{GW:periods}. 
To begin with, we fix an ideal triangulation $\calT$ with $N$ tetrahedra of
a 3-manifold $M$ whose torus boundary is marked by a pair of a meridian and
longitude. The building block of the 3D-index is the tetrahedron index
$\ID(m,e)(q) \in \BZ[[q^{1/2}]]$ defined by
\be
\label{ID}
\ID(m,e)(q)=\sum_{n=(-e)_+}^\infty (-1)^n \frac{q^{\frac{1}{2}n(n+1)
-\left(n+\frac{1}{2}e\right)m}}{(q;q)_n(q;q)_{n+e}}
\ee
where $e_+=\max\{0,e\}$ and $(q;q)_n=\prod_{i=1}^n (1-q^i)$. If we wish, we can
sum in the above equation over the integers, with the understanding that
$1/(q;q)_n=0$ for $n < 0$.

The rotated 3D-index is given by
\be
\label{Irotdef}
\Irot_\calT(n,n')(q) = \sum_{k \in \BZ^N} S_\calT(k,n,n')(q)
\ee
where 
\begin{small}
\be
\label{Sdef}
S_\calT(k,n,n')(q) =
(-q^{1/2})^{\mb \nu \cdot k-(n-n')\nu_\lambda} q^{k_N(n+n')/2}
\prod_{j=1}^N \ID(\lambda''_j(n-n') -b_j \cdot k, -\lambda_j(n-n') + a_j \cdot k)(q) 
\ee
\end{small}
is assembled out of a product of tetrahedra indicies $\ID$ 
evaluated to linear forms that depend on the Neumann--Zagier matrices
$(A|B)$ of $\calT$. The detailed definition of the Neumann--Zagier matrices is given
in Appendix~\ref{app.indexdef}.

Note that the degree $\delta(\ID(m,e))$ of the tetrahedron index
is a nonnegative piecewise quadratic function of $(m,e)$
\be
\label{degID}
\delta(\ID(m,e)) = \frac{1}{2}\left( m_+(m+e)_+ + (-m)_+ e_+ + (-e)_+(-e-m)_+
+\max\{0,m,-e\} \right) \,.
\ee
It follows that for 1-efficient triangulations (see~\cite{GHRS}) the degree of
the summand in~\eqref{Irotdef} is bounded below by a positive constant times
$\max\{|k_1|,|k_2|,\dots,|k_N|\}$,
thus the sum in~\eqref{Irotdef} is a well-defined element of $\BZ((q^{1/2}))$. 

The topological invariance of the 3D-index is a bit subtle, since the definition 
requires 1-efficient ideal triangulations,
and the latter are not known to be connected under 2--3 Pachner moves. Nonetheless,
in~\cite{GHRS}, it was shown that the 3D-index (and likewise, its rotated version) is
a topological invariant of cusped hyperbolic 3-manifolds. An alternative proof of
this fact was given in~\cite{GK:mero}, where the rotated 3D-index was reformulated
in terms of a meromorphic function of two variables.

\subsection{Factorization and holomorphic blocks}
\label{sub.factorization}

From its very definition as a sum of proper $q$-hypergeometric series, it follows
that $\Irot_\calT(n,n')(q)$ is a $q$-holonomic function of $n$ and $n'$~\cite{WZ, AB}.
But more is true. The rotated 3D-index factorizes into a sum of a product of pairs
of colored holomorphic blocks. This holomorphic block factorization is a well-known
phenomenon explained in~\cite{Beem}, and most recently in~\cite{GW:periods} whose
presentation we will follow. Let us recall how this works. We can assemble the
collection $\Irot_\calT(n,n')(q)$ of $q$-series indexed by pairs of integers into a
$\BZ\times\BZ$ matrix $\Irot_\calT(q)$ whose $(n,n')$ entry is $\Irot_\calT(n,n')(q)$.
Then, in~\cite{GW:periods} we explained the origin of the following conjecture for
the rotated 3D-index. 

\begin{conjecture}
\label{conj.0}
For every $1$-efficient triangulation $\calT$ there exists a palindromic linear
$q$-difference operator $\Ahat_\calT$ of order $r$ with a fundamental solution
$\BZ \times r$ matrix $H_\calT(q)$ and a symmetric, invertible $r \times r$ matrix
$B_\calT$ with rational entries such that
\be
\label{IrotH}
\Irot_\calT(q) = H_\calT(q) B_\calT H_\calT(q^{-1})^t \,.
\ee
\end{conjecture}
When the triangulation is fixed and clear, we will drop it from the notation.
If we denote the $(n,\a)$ entry of $H_\calT(q)$ whose $(n,\a)$ entry by $h^{(\a)}_n(q)$,
these functions are the so-called colored holomorphic blocks. It follows that
the matrix $H(q)$ is a (properly normalized)
fundamental solution to a pair of $q$-difference equations
\be
\label{AH}
\Ahat_\calT(M_+,L_+) H(q) = 0, \qquad \Ahat_\calT(M_-,L_-) H(q^{-1}) =
0,
\ee
where the operators act respectively by
\begin{equation}
  \begin{gathered}
  M_+ h_n^{(\alpha)}(q) =q^n h_n^{(\alpha)}(q) ,\quad L_+
  h_n^{(\alpha)}(q) = h_{n+1}^{(\alpha)}(q)\\
  M_- h_n^{(\alpha)}(q^{-1}) =q^{-n }h_n^{(\alpha)}(q^{-1}) ,\quad L_-
  h_n^{(\alpha)}(q^{-1}) = h_{n+1}^{(\alpha)}(q^{-1}).
  \end{gathered}
\end{equation}
Consequently the rotated 3D-index satisfies a pair of (left and right)
linear $q$-difference equations \be
\label{AhatIrot}
\Ahat_\calT(M_+,L_+) \Irot_\calT = \Ahat_\calT(M_-,L_-) \Irot_\calT =0 
\ee
acting in a decoupled way on each of the rows and columns of $\Irot_\calT$.

The factorization~\eqref{IrotH} of the rotated 3D-index and the left and right linear
$q$-difference equations~\eqref{AhatIrot} imply the following.

\begin{corollary}
\label{cor.Irot}(of Conjecture~\ref{conj.0})
The rotated 3D-index $\Irot_\calT(q)$ is uniquely determined by
\begin{itemize}
\item[(1)]
  the $r \times r$ matrix $\Irot_\calT(q)[r]$ and
\item[(2)] the pair of linear $q$-difference equations~\eqref{AhatIrot}.
\end{itemize}
\end{corollary}

\noindent
Here, $\Irot_\calT(q)[r]$ denotes the $r \times r$ matrix
$(\Irot_\calT(n,n')(q))$ for $0 \leq n,n' \leq r-1$. 

The holomorphic blocks satisfy the symmetry
\be
\label{hansym}
h^{(\a)}_{\calT,n}(q) = h^{(\a)}_{\calT,-n}(q) 
\ee
for all $\a$ and all integers $n$, which together with Equation~\eqref{IrotH}
implies the symmetries
\be
\label{Isym}
\Irot_\calT(n,n')(q)=\Irot_\calT(n,-n')(q)=\Irot_\calT(-n,n')(q)
=\Irot_\calT(-n,-n')(q) \,,
\ee
and 
\be
\label{Ikmsym}
\Irot_\calT(n,n')(q^{-1}) = \Irot_\calT(n',n)(q) \,,
\ee
for the rotated 3D-index.

Let us finally mention that the colored holomorphic blocks can be computed by the
limit as $x \to 1$
\be
\label{IrotB}
\Irot_\calT(n,n')(q) =
\lim_{x\rightarrow1}\sum_\a B^{(\a)}_\calT(q^{-n'}x^{-1};q^{-1})
B^{(\a)}_\calT(q^{n}x;q)\,.
\ee
of the $x$-deformed holomorphic blocks $B^{(\a)}_\calT(x;q)$ and the latter can be
determined from a factorization of an appropriate state-integral.

\subsection{Descendants}
\label{sub.desc}
  
There is an important $\BQ(q)$-linear action of $\BW_q(\calT)$ on the set of
functions $S_\calT(k,n,n')(q)$ giving rise to a map
\be
\label{IrotW}
\BW_q(\calT) \to \BZ((q^{1/2}))^{\BZ^N \times\BZ^2}
\ee
which descends to a push-forward $\BQ(q^{1/2})$-linear map
\be
\label{IrotO}
\calM(\calT) \to \BZ((q^{1/2}))^{\BZ^2}, \qquad \calO \mapsto \Irot_{\calT,\calO} \,.
\ee
Concretely, when $\calO=\prod_{j=1}^N \hz_j^{\a_j} (\hz_j^{''})^{\b_j}$, we have
\be
\label{IrotOdef}
\Irot_{\calT,\calO}(n,n')(q) =
\sum_{k \in \BZ^N} (\calO \circ S_\calT)(k,n,n')(q) \,,
\ee
where
\begin{small}
\be
\label{SOdef}
\begin{aligned}
  (\calO \circ S_\calT)(k,n,n')(q) & =
(-q^{1/2})^{\mb \nu \cdot k-(n-n')\nu_\lambda} q^{k_N(n+n')/2+L_\calO(n,n',k)}  \\
  & \hspace{-2.5cm}
\times \prod_{j=1}^N 
\ID(\lambda''_j(n-n') -b_j \cdot k + \beta_j,
-\lambda_j(n-n') + a_j \cdot k - \alpha_j)(q) \,,
\end{aligned}
\ee
\be
\label{Lnk}
L_\calO(n,n',k) = \frac{1}{2}\sum_{j=1}^N
\big(\alpha_j(\lambda''_j n-\lambda''_j n' -b_j \cdot k)
+ \beta_j(-\lambda_j n + \lambda_j n' + a_j \cdot k) - \alpha_j\beta_j\big)
\ee
\end{small}

\noindent
This action was written down explicitly in~\cite[Eqn.(104)]{Agarwal}.
The symmetries of the tetrahedron index~\cite[Eqns.(136)]{DGG1} imply that
the three Lagrangian operators given in Equation~\eqref{3lag} annihilate
$S_\calT(k,n,n')(q)$, and thus the sum $\Irot_\calT(n,n')(q)$. In addition, the
insertion $\calE_i$ corresponding to the $i$-th edge (for $i=1,\dots, N-1$) when quantized as in \cite{Dimofte:riemann} satisfies
\be
(\calE_i \circ S_\calT)(k,n,n')=
q S_\calT(k- e_i,n,n')
\ee
Summing over $k$, this implies that
$\calE_i- q$ annihilates 
$\Irot_\calT(n,n')(q)$. Thus, $\Irot_{\calT,\calO}(n,n')(q)$ is well-defined
for all $\calO \in \calM(\calT)$, justifying the strange quotient given in
Equation~\eqref{strangequotient}. 
Note that the action of the edge operators considered in~\cite{Dimofte:riemann} differs
by factor of $q$ from that of~\cite[Eqn.(130)]{Agarwal}. 

Our conjecture relates the colored holomorphic blocks and the
rotated 3D-index of $\calT$ to those of $(\calT,\calO)$. Simply put, it asserts that
inserting $\calO$ simply changes the invariants ($\BZ((q^{1/2}))$-series) by
multiplication of a matrix of rational functions, and changes the left $q$-difference
equation whereas it preserves the right one. This implies that the $\BQ(q^{1/2})$-span
of the collection $\{\Irot_{\calT,\calO}(q) \, | \, \calO \in \calM(\calT)\}$ is
a finite dimensional $\BQ(q^{1/2})$-vector space.

Fix a 1-efficient ideal triangulation $\calT$ of a 1-cusped 3-manifold $M$.

\begin{conjecture}
\label{conj.1}
For every $\calO \in \calM(\calT)$
\begin{itemize}
\item[(a)]
  there exists a linear $q$-difference operator $\Ahat_{\calT,\calO}$ with a
  fundamental solution matrix $H_{\calT,\calO}(q)$ such that 
\be
\label{IrotHO}
\Irot_{\calT,\calO}(q) = H_{\calT,\calO}(q) B_\calT H_\calT(q^{-1})^t \,,
\ee
\item[(b)]
  there exists $Q_{\calT,\calO}(q) \in \mathrm{GL}_r(\BQ(q^{1/2}))$ such that
\be
\label{IHQ}
\Irot_{\calT,\calO}[r] = Q_{\calT,\calO} \Irot[r], \qquad
H_{\calT,\calO}[r] = Q_{\calT,\calO} H[r] \,.
\ee
\end{itemize}
\end{conjecture}

The above conjecture implies the following.

\begin{corollary}
\label{cor.IrotO}(of Conjecture~\ref{conj.1})
The rotated 3D-index $\Irot_{\calT,\calO}(q)$ is uniquely determined by
\begin{itemize}
\item[(1)]
  the $r \times r$ matrices $\Irot_\calT(q)[r]$ and $Q_{\calT,\calO}(q)$
\item[(2)] the pair of linear $q$-difference equations $\Ahat_{\calT,\calO}$ and
  $\Ahat_\calT$.
\end{itemize}
\end{corollary}

Another corollary of the above conjecture concerns the descendants of the
rotated 3D-index, analogous to the descendants of the colored Jones polynomial of
a knot defined in~\cite{GK:desc} and the descendants of the holomorphic blocks
defined in~\cite[Eqn.(13), App.A]{GGM}. To phrase it, let
\be
\label{Dspan}
D\Irot_\calT =
\text{Span}_{\BQ(q^{1/2})} \{ \Irot_\calT(n,n')(q) \, | \, n,n' \in \BZ \}
\ee
denote the $\BQ(q^{1/2})$-span of the elements $\Irot_\calT(n,n')$ of the ring
$\BQ((q^{1/2}))$. Note that $D\Irot_\calT$ is a finite dimensional vector space
over the field $\BQ(q^{1/2})$. Likewise, one defines $\Irot_{\calT,\calO}$. 
The next corollary justifies the title of the paper. 

\begin{corollary}
\label{cor.Irotdesc}(of Conjecture~\ref{conj.1})
We have:
\be
\label{IOspan}
\cup_{\calO \in \calM(\calT)} D\Irot_{\calT,\calO} = D\Irot_\calT \,. 
\ee
\end{corollary}

In other words, the descendants $D\Irot_{\calT,\calO}$ of the rotated 3D-index 
$D\Irot_{\calT}$ are expressed effectively by a finite-size matrix with entries
in $\BQ(q^{1/2})$. 

We now formulate a relative version of the AJ-Conjecture. 
Let $\Ahat(M,L)|_{q=1}=A(M,L)$ denote the classical limit of a linear $q$-difference
equation. The AJ-Conjecture~\cite{Ga:AJ} relates the classical limit of the
$\Ahat$-polynomial with the $A$-polynomial of a knot given in~\cite{CCGLS}.

\begin{conjecture}
\label{conj.2}
For every $\calO \in \calM(\calT)$, we have
\be
A_{\calT,\calO}(M,L) =_M A_\calT(M,L) 
\ee
where $=_M$ means equality up to multiplication by a nonzero function of $M$. 
\end{conjecture}

\subsection{Asymptotics}
\label{sub.asy}

A consequence of Conjecture~\eqref{conj.1} (and Equation~\eqref{IrotHO}) is that
the all-order asymptotics of the colored holomorphic blocks $h^{(\a)}_{\calT,\calO,n}(q)$
and the $\Irot_{\calT,\calO}(n,n')(q)$ are a $\BQ(q)$-linear combination of those of
$h^{(\a)}_{\calT,n}(q)$ and $\Irot_\calT(n,n')(q)$, respectively. The asymptotics
of the latter were studied in detail in~\cite{GW:periods}. A corollary of this 
and Conjecture~\ref{conj.2} is a resolution and an explanation from first
principles, of the quantum length conjecture of~\cite{Agarwal}. 


\section{Examples}
\label{sec.examples}


In this section we illustrate our conjectures with the case of the three
simplest  hyperbolic knots, the $4_1$ (figure eight) knot, the $5_2$ knot
and the $(-2,3,7)$ pretzel knot. 

\subsection{The $4_1$ knot and its rotated 3D-index}
\label{sub.Irot41}

The complement of the $4_1$ knot has an ideal triangulation with two tetrahedra.
Using the gluing equation matrices 
\be
\mb G=\begin{pmatrix} 
 2 & 2 \\
 0 & 0 \\
 1 & 0 \\
 1 & 1 \\
  \end{pmatrix}
  , \qquad
\mb G'=\begin{pmatrix} 1 & 1 \\
 1 & 1 \\
 0 & 0 \\
 1 & -1 \\ 
  \end{pmatrix}
  , \qquad
\mb G''=\begin{pmatrix}  0 & 0 \\
 2 & 2 \\
 0 & -1 \\
 1 & -3 \\ 
  \end{pmatrix}
, \qquad  
\ee
with the conventions explained in Appendix~\ref{app.indexdef}, we obtain the matrices 
\be
A = \begin{pmatrix} 1 & 1 \\ 
  1 & 0  
\end{pmatrix}, B = \begin{pmatrix} -1 & -1 \\ 
  0 & -1  
\end{pmatrix}, \nu = \begin{pmatrix} 0 \\ 
  0  
\end{pmatrix}
\ee
in terms of which, the rotated 3D-index is given by
\be
\label{Irot41}
\Irot_{4_1}(n,n')(q) = \sum_{k_1,k_2 \in \BZ}
q^{k_2(n+n^{\prime})/2}\ID(k_1,k_1+k_2)(q) \ID(k_1+k_2-n+n^\prime,k_1-n+n^\prime)(q)
\ee
where $\ID$ is the tetrahedron index given in~\eqref{ID}. (The above formula
agrees with~\cite[Eqn.(108)]{Agarwal} after a shift $k_1 \mapsto k_1 - k_2$).  
Using Equation~\eqref{degID}, it follows that the degree of the summand
in~\eqref{Irot41} is bounded below by a positive constant times
$\max\{|k_1|,|k_2|\}$, thus the sum in~\eqref{Irot41} is a well-defined element
of $\BZ((q^{1/2}))$. 

\subsection{Factorization}
\label{sub.41GW}

In this section we briefly summarize the properties of the rotated 3D-index of
the $4_1$ knot following~\cite{GW:periods}, namely its factorization in terms of
colored holomorphic blocks, the linear $q$-difference equation, their symmetries
and their asymptotics. All the functions in this section involve the knot $4_1$,
which we suppress from the notation.

The rotated 3D-index is given by~\cite[Prop.9]{GW:periods}
  \be
  \label{41.irot}
  \Irot_{4_1}(n,n')(q)
  =
  -\frac{1}{2}h_{4_1,n'}^{(1)}(q^{-1})h_{4_1,n}^{(0)}(q)
  +\frac{1}{2}h_{4_1,n'}^{(0)}(q^{-1})h_{4_1,n}^{(1)}(q)\, \qquad (n , n' \in \BZ)
\ee
with the colored holomorphic blocks $h^{(0)}_{4_1,n}(q)$ and $h^{(1)}_{4_1n}(q)$ given
in the Appendix~\ref{app.41}. 

The colored holomorphic blocks satisfy the symmetries
\be
\label{41hsym1}
h_{4_1,n}^{(0)}(q^{-1})=h_{4_1,n}^{(0)}(q), \qquad
h_{4_1,n}^{(1)}(q^{-1})=-h_{4_1,n}^{(1)}(q) \,,
\ee
and
\be
\label{41hsym2}
h_{4_1,-n}^{(\a)}(q) = h_{4_1,n}^{(\a)}(q), \qquad \a=0,1 \,,
\ee
and the linear $q$-difference equation~\cite[Eqn.(63)]{GW:periods}
\be
\label{41.Ahath}
P_{4_1,0}(q^n,q)h^{(\a)}_n(q)+P_{4_1,1}(q^n,q)h^{(\a)}_{n+1}(q)
+P_{4_1,2}(q^n,q)h^{(\a)}_{n+2}(q)=0 \qquad (\a=0,1, \, \, n \in \BZ) 
\ee
where
\be
\label{Ahat41}
\begin{aligned}
  P_{4_1,0}(x,q)&=q^{2}x^2(q^{3}x^2-1)\,,\\
  P_{4_1,1}(x,q)&=-q^{1/2}(1-q^{2}x^2)(1-qx-qx^2-q^{3}x^2-q^{3}x^3+q^{4}x^4)\,,\\
  P_{4_1,2}(x,q)&=q^{3}x^2(-1+qx^2) \,.
\end{aligned}
\ee
We denote the corresponding operator of the $q$-difference equation~\eqref{41.Ahath}
by $\Ahat_{4_1}(x,\s,q)=\sum_{j=0}^2 P_{4_1,j}(x,q)\s^j$.

\subsection{Defects}
\label{sub.Irot41insertions}

We now consider two defects. The first one is the element
\be
\label{O1}
\calO=-\wh y^{-1} -\wh z^{-1} + \wh y^{-1} \wh z^{-1} \in \calM(\calT) 
\ee
from~\cite[Eqn.(81)]{Agarwal}. Computing the values of
$\Irot_{4_1}(n,n')(q)$ and $\Irot_{4_1,\calO}(n,n')(q)$ for
$0 \leq n,n' \leq 1$ up to $O(q^{121})$, we find out that the
$2 \times 2$ matrices
\begin{tiny}
  \begin{multline*}
  \Irot_{4_1}(q)[2]= \\
  \begin{pmatrix}
1 - 8 q - 9 q^2 + 18 q^3 + 46 q^4 + 90 q^5 + 62 q^6 + 
  10 q^7 + \dots & -q^{-1/2} + q^{1/2} - q^{3/2} + 6 q^{5/2} + 20 q^{7/2} + 
  29 q^{9/2} + 25 q^{11/2} 
  + \dots
  \\
  -q^{-1/2} + q^{1/2} - q^{3/2} + 6 q^{5/2} + 
  20 q^{7/2} + 29 q^{9/2} 
  + \dots &
 2 q + 2 q^2 + 7 q^3 + 8 q^4 + 3 q^5 - 22 q^6 - 67 q^7 + \dots
\end{pmatrix}
\end{multline*}
\end{tiny}
and
\begin{tiny}
\begin{multline*}
\Irot_{4_1,\calO}(q)[2]= \\
\begin{pmatrix}
-3 + 15 q + 24 q^2 - 15 q^3 - 69 q^4 - 174 q^5 - 183 q^6 - 165 q^7
+ \dots &
 2 q^{-1/2} - q^{1/2} + 4 q^{3/2} - 7 q^{5/2} - 34 q^{7/2} - 
 64 q^{9/2} 
 + \dots
 \\
 q^{-3/2} - q^{-1/2} - q^{1/2} + q^{3/2} - 5 q^{5/2} - 26 q^{7/2} - 
 48 q^{9/2} 
 + \dots
 & -1 - 2 q - 4 q^2 - 9 q^3 - 17 q^4 - 13 q^5 + 10 q^6 + 77 q^7
 + \dots
\end{pmatrix}
\end{multline*}
\end{tiny}
satisfy
\be
\label{cancel}
(q-1)\Irot_{4_1,\calO}(q)[2](\Irot(q)_{4_1}[2])^{-1} =
\begin{pmatrix} 2-q & -q^{1/2} \\ 
  q^{1/2} & -q-1+q^{-1}  
\end{pmatrix} + O(q^{121}) 
\ee
illustrating the dramatic collapse of the $q$-series into short rational
functions of $q^{1/2}$. This implies that the matrix $Q_{4_1,\calO}(q)$ is given by
\be
\label{Q1}
Q_{4_1,\calO}(q) = \frac{1}{q-1}
\begin{pmatrix} 2-q & -q^{1/2} \\ 
  q^{1/2} & -q-1+q^{-1}  
\end{pmatrix} 
\ee
with $\det(Q_{4_1,\calO})(q)=1+2 q^{-1}$.

After computing the values of $\Irot_{4_1,\calO}(n,0)(q) + O(q^{120})$ for $n=0,\dots,10$
and finding a short linear recursion among three consecutive values, and further
interpolating for all $n$, we found out that the left $\Ahat$-polynomial 
of $\Irot_{4_1,\calO}(q)$ is given by 
$\Ahat_{4_1,\calO}(x,\s,q)=\sum_{j=0}^2 P_{4_1,\calO,j}(x,q)\s^j$ where
\be
\label{Ahat41O}
\begin{aligned}
  P_{4_1,\calO,0}(x,q)&=q^{3/2} x^2 (-1 + q^3 x^2) (1 + q x + q^3 x^2)\,,\\
  P_{4_1,\calO,1}(x,q)&=(-1 + q x) (1 + q x) \\
  & (1 + x - q x - q x^2 - q^3 x^2 - q x^3 - 
   2 q^3 x^3 - q^5 x^3 - q^3 x^4 - q^5 x^4 + q^4 x^5 - q^5 x^5 + 
   q^6 x^6)\,,\\
  P_{4_1,\calO,2}(x,q)&=q^{7/2} x^2 (-1 + q x^2) (1 + x + q x^2)\,.
\end{aligned}
\ee
The $\Ahat_{4_1,\calO}$ polynomial is palindromic, and together with the skew-symmetry
of the $Q_{4_1,\calO}(q)$ matrix, it follows that the colored holomorphic blocks
$h^{(0)}_{4_1,\calO,n}(q)$ and $h^{(1)}_{4_1,\calO,n}(q)$ satisfy the
symmetries~\eqref{41hsym1} and~\eqref{41hsym2}.

When we set $q=1$, we obtain
\be
\Ahat_{4_1,\calO}(x,\s,1) = 2(x^2-1)(x^2+x+1) \Ahat_{4_1}(x,\s,1)
\ee
confirming Conjecture~\ref{conj.2}.

Equation~\eqref{Q1} and the recursion~\eqref{Ahat41O} imply that for all integers
$n$ and $n'$, $\Irot_{4_1,\calO}(n,n')(q)$ is a $\BQ(q^{1/2})$-linear combination of
the three $q$-series $\Irot_{4_1}(0,0)(q)$, $\Irot_{4_1}(0,1)(q)$ and
$\Irot_{4_1}(1,0)(q)$. For instance, Equation~\eqref{IHQ} implies that
\be
\label{ex41}
\Irot_{4_1,\calO}(0,0)(q) = \tfrac{1}{q-1}
((2-q)\Irot_{4_1}(0,0)(q) -q^{\tfrac{1}{2}} \Irot_{4_1}(0,1)(q)) 
\ee
and likewise for other values of $\Irot_{4_1,\calO}(n,n')(q)$. This reduces the
problem of the asymptotic expansion of $\Irot_{4_1,\calO}(n,n')(q)$ for
$q=e^{2\pi i \tau}$ to all orders in $\tau$ as $\tau$ tends to zero in a ray (nearly
vertically, horizontally, or otherwise) to the problem of the asymptotics of
colored holomorphic blocks and of the rotated 3D-index. This problem was studied
in detail and solved in the work of Wheeler and the second
author~\cite[Sec.5.7,5.8]{GW:periods} for the $4_1$ knot.

As a second example, consider the element 
\be
\label{O2}
\calO_2=\wh y^{-1} \in \calM(\calT) \,.
\ee
Repeating the above computations, we find out that the matrix $Q_{4_1,\calO_2}(q)$
is given by
\be
\label{Q2}
Q_{4_1,\calO_2}(q) = \frac{1}{q-1}
\begin{pmatrix} -1 & q^{1/2} \\ 
  -q^{1/2} & -q^2+2 q + 1 -q^{-1}  
\end{pmatrix} 
\ee
with $\det(Q_{4_1,\calO_2})(q)=1+q^{-1}$, and that the left $\Ahat$-polynomial 
of $\Irot_{4_1,\calO_2}(q)$ is given by 
$\Ahat_{4_1,\calO_2}(x,\s,q)=\sum_{j=0}^2 P_{4_1,\calO_2,j}(x,q)\s^j$ where
\be
\label{Ahat41O2}
\begin{aligned}
P_{4_1,\calO_2,0}(x,q)&=q^{3/2} x^2 (-1 + q^2 x) (1 + q^2 x)
\,, \\
P_{4_1,\calO_2,1}(x,q)&=(-1 + q^3 x^2) (1 - q x - q^2 x^2 - q^4 x^2 - q^4 x^3 + q^6 x^4)
\,, \\
P_{4_1,\calO_2,2}(x,q)&=q^{7/2} x^2 (-1 + q x) (1 + q x)
\,.
\end{aligned}
\ee
In this case, we lose the Weyl-invariance symmetry of the colored holomorphic blocks,
but we retain the AJ Conjecture~\ref{conj.2} since
\be
\Ahat_{4_1,\calO_2}(x,\s,1) = (x^2-1) \Ahat_{4_1}(x,\s,1) \,.
\ee

\subsection{The $5_2$ knot and its rotated 3D-index}
\label{sub.Irot52}

The complement of the $5_2$ knot has an ideal triangulation with three tetrahedra.
Using the gluing equation matrices 
\be
\mb G=\begin{pmatrix} 1 & 1 & 1 \\
 0 & 0 & 0 \\
 1 & 1 & 1 \\
 -1 & 0 & 0 \\
 3 & 2 & 1 \\
  \end{pmatrix}
  , \qquad
\mb G'=\begin{pmatrix} 0 & 2 & 0 \\
 1 & 0 & 1 \\
 1 & 0 & 1 \\
 0 & 0 & 0 \\
 1 & 2 & 1 \\
  \end{pmatrix}
  , \qquad
\mb G''=\begin{pmatrix} 1 & 0 & 1 \\
 1 & 2 & 1 \\
 0 & 0 & 0 \\
 0 & 1 & 0 \\
 -1 & 0 & 3 \\
  \end{pmatrix}
, \qquad  
\ee
with the conventions explained in Appendix~\ref{app.indexdef}, we obtain the matrices 
\be
A = \begin{pmatrix} 1 & -1 &1 \\ 
  -1 & 0 & -1\\
  -1 & 0 & 0
\end{pmatrix},\quad B = \begin{pmatrix} 1 & -2 &1 \\ 
  0 & 2 & 0\\
  0 & 1 & 0
\end{pmatrix},\quad \nu = \begin{pmatrix} 0 \\ 
  0 \\
  0
\end{pmatrix}\,.
\ee
The rotated 3D-index is given by 
\begin{small}
\be
\label{Irot52}
\begin{aligned}
\Irot_{5_2}(n,n')(q) =
&\sum_{k_1,k_2,k_3 \in \BZ} q^{k_3(n+n^\prime)/2}
\ID(k_1 - k_2,k_3+k_2+n-n^\prime)\\
& \hspace{-2cm} \times \ID(-k_1+2k_2-n+n^\prime, k_3 +2k_1 - 2k_2+ n - n^\prime)
\ID(k_3 + k_1 -k_2 +n -n^\prime,k_2 - 2n+2n^\prime)\,.
\end{aligned}
\ee
\end{small}
Equation~\eqref{degID} implies that the degree of the summand
in~\eqref{Irot52} is bounded below by a positive constant times
$\max\{|k_1|,|k_2|,|k_3|\}$, thus the sum in~\eqref{Irot52} is a well-defined element
of $\BZ((q^{1/2}))$. 

\subsection{Factorization}
\label{sub.52GW}

The $5_2$ knot has three colored holomorphic blocks $h^{(\a)}_n(q)$ for $\a=0,1,2$,
$n$ an integer and $q$ a complex number $|q| \neq 1$, whose definition in terms of
$q$-hypergeometric series was given in~\cite[App.A]{GW:periods} and reproduced for
the convenience of the reader in Appendix \ref{app.52}. 
The rotated 3D-index is given by~\cite[Prop.13]{GW:periods}
\be
\label{52.irot}
\Irot_{5_2}(n,n')(q)
=
-\frac{1}{2}h_{5_2,n'}^{(0)}(q^{-1})h_{5_2,n}^{(2)}(q)
-h_{5_2,n'}^{(1)}(q^{-1})h_{5_2,n}^{(1)}(q)
-\frac{1}{2}h_{5_2,n'}^{(2)}(q^{-1})h_{5_2,n}^{(0)}(q)\,.
\ee

The colored holomorphic blocks satisfy the symmetries
\be
\label{52hsym2}
h_{5_2,-n}^{(\a)}(q) = h_{5_2,n}^{(\a)}(q), \qquad \a=0,1,2 \,.
\ee
and the linear $q$-difference equation~\cite[Eqn.(63)]{GW:periods}
\be
\label{52.Ahat.colored}
\begin{aligned}
P_{5_2,0}(q^n,q)h^{(\alpha)}_{n}(q)
+P_{5_2,1}(q^n,q)h^{(\alpha)}_{n-1}(q)
+P_{5_2,2}(q^n,q)h^{(\alpha)}_{n-2}(q)
+P_{5_2,3}(q^n,q)h^{(\alpha)}_{n-3}(q)&=0\,,
\end{aligned}
\ee
for all $\a=0,1,2$ and all integers $n$, where~\cite[Eqn.(126)]{GW:periods}
\be
\label{Ahat52coeffs}
\begin{tiny}
\begin{aligned}
P_{5_2,0}(x,q)&=-q^{-2}x^2(1-q^{-2}x)(1+q^{-2}x)(1-q^{-5}x^2)\,,\\
P_{5_2,1}(x,q)&=q^{3/2}x^{-3}(1-q^{-1}x)(1+q^{-1}x)(1-q^{-5}x^2) \\
& \quad \cdot
(1-q^{-1}x-q^{-1}x^2-q^{-4}x^2+q^{-2}x^2+q^{-3}x^2+q^{-2}x^3
+q^{-5}x^3+q^{-5}x^4+q^{-5}x^4-q^{-6}x^5)\,,\\
P_{5_2,2}(x,q)&=q^{5}x^{-5}(1-q^{-2}x)(1+q^{-2}x)(1-q^{-1}x^2) \\
& \quad \cdot
(1-q^{-2}x-q^{-2}x-q^{-2}x^2-q^{-5}x^2+q^{-4}x^3+q^{-7}x^3-q^{-5}x^3
-q^{-6}x^3+q^{-7}x^4-q^{-9}x^5)\,,\\
P_{5_2,3}(x,q)&=q^{\frac{11}{2}}x^{-5}(1-q^{-1}x)(1+q^{-1}x)(1-q^{-1}x^2)\,.
\end{aligned}
\end{tiny}
\ee

\subsection{Defects}
\label{sub.Irot52insertions}

We now consider two defects $\calO_1$ and $\calO_2$ given by
\be
\label{52O}
\begin{aligned}
\calO_1 &=\wh z_1 \\
\calO_2 &=\wh z_1 + \wh z_3 \,.
\end{aligned}
\ee
Computing the $3 \times 3$ matrix of the rotated 3D-index with and without
insertion up to $O(q^{81})$, and dividing one matrix by another, we found out
that the corresponding $3 \times 3$ matrices $Q_{\calO_j}(q)+ O(q^{81})$ for
$j=1,2$ are given by
\begin{tiny}
\label{52cancel1}
\begin{multline}
  \Irot_{5_2,\calO_1}(q)[3](\Irot_{5_2}(q)[3])^{-1} = \frac{1}{(1-q^2)(1-q^3)} \cdot
\\
\begin{pmatrix}
  -q^2 - q^3 - q^4 & 
  q^{1/2} - q^{3/2} + q^{7/2} + 2 q^{9/2} + 2 q^{11/2} - q^{13/2} & -q^7 \\
  -q^{3/2} - q^{5/2} - q^{7/2} & 
  1 - q + q^3 + 2 q^4 + 2 q^5 - q^6 & -q^{13/2} \\
  -1 - q^{-2} - q^{-1} &
  -q^{-5/2} + 2 q^{-3/2} + 2 q^{-1/2} + q^{1/2} - q^{5/2} + q^{7/2} & -q^3
\end{pmatrix} + O(q^{81}) 
\end{multline}
\end{tiny}
and
\begin{tiny}
\label{52cancel2}
\begin{multline}
\Irot_{5_2,\calO_2}(q)[3](\Irot_{5_2}(q)[3])^{-1} = \frac{1}{(1-q^2)(1-q^3)} \cdot
\\
\begin{pmatrix}
-q - 2 q^2 - q^3 + q^5 &
  q^{1/2} + q^{5/2} + q^{7/2} + q^{9/2} + q^{11/2} - q^{13/2} & -q^7
\\ -q^{-1/2} - q^{3/2} - q^{7/2} &
  4 - q^{-1} - q - q^2 - q^3 + q^4 + 5 q^5 - 2 q^6 & 
  q^{11/2} - 2 q^{13/2}
\\ -2 + q^{-4} + q^{-3} - q^{-2} - 2 q^{-1} &
  q^{-9/2} - 2 q^{-7/2} - 4 q^{-5/2} + 2 q^{-3/2} + 4 q^{-1/2} + 
   4 q^{1/2} - q^{3/2} - 2 q^{5/2} + 2 q^{7/2} & 1 + q - q^2 - 2 q^3
\end{pmatrix} + O(q^{81}) 
\end{multline}
\end{tiny}
illustrating Corollary~\ref{cor.Irotdesc} of Conjecture~\ref{conj.1}.

\subsection{The $(-2,3,7)$-pretzel knot}
\label{sub.237}

As a final experiment, we studied the rotated 3D-index of the $(-2,3,7)$ pretzel-knot.
This knot is interesting in several ways, and exhibits behavior of general hyperbolic
knots. The complement of the $(-2,3,7)$-pretzel knot is geometrically similar
to that of the $5_2$ knot, i.e., both are obtained by the gluing of three 
three ideal tetrahedra, only put together in a combinatorially different way. Thus,
the $5_2$ and $(-2,3,7)$ pretzel knots have the same cubic trace field, and the same
real volume. But the similarities end there. The $5_2$ knots has three boundary
parabolic $\PSL_2(\BC)$-representations, all Galois conjugate to the geometric one.
On the other hand, one knows from~\cite{GZ:kashaev} and~\cite{GZ:qseries} that the
$(-2,3,7)$-pretzel knot has 6
colored holomorphic blocks, corresponding to the fact that the $(-2,3,7)$-pretzel
knot has 6 boundary parabolic representations, three coming from the Galois orbit of
the geometric $\PSL_2(\BC)$-representation (defined over the cubic trace field of 
discriminant $-23$) and three more coming from the Galois orbit of a
$\PSL_2(\BC)$-representation defined over the totally real abelian
field $\BQ(\cos(2 \pi/7))$. Although~\cite{GZ:qseries} gives explicit expressions
for the $6 \times 6$ matrices of the holomorphic blocks (inside and outside the
unit disk), the colored holomorphic blocks have not been computed, partly due to the
complexity of the calculation.

Going back to the 3D-index of the $(-2,3,7)$ knot, the gluing equation matrices
are 
\be
\mb G=\begin{pmatrix} 1 & 1 & 1 \\
 1 & 0 & 0 \\
 0 & 1 & 1 \\
 0 & 0 & -1 \\
 -1 & 1 & -18 \\
  \end{pmatrix}
  , \qquad
\mb G'=\begin{pmatrix}1 & 0 & 0 \\
 0 & 2 & 2 \\
 1 & 0 & 0 \\
 0 & 0 & 0 \\
 1 & -1 & -2 \\
  \end{pmatrix}
  , \qquad
\mb G''=\begin{pmatrix}  0 & 1 & 0 \\
 2 & 1 & 0 \\
 0 & 0 & 2 \\
 2 & 0 & 0 \\
 35 & 1 & 0 \\
  \end{pmatrix}
, \qquad  
\ee
with the conventions explained in Appendix~\ref{app.indexdef}, from which we
obtain that
\be
A = \begin{pmatrix} 0 & 1 &1 \\ 
  1 & -2 & -2\\
  0 & 0 & -1
\end{pmatrix},\quad B = \begin{pmatrix} -1 & 1 &0 \\ 
  2 & -1 & -2\\
  2 & 0 & 0
\end{pmatrix},\quad \nu = \begin{pmatrix} 1 \\ 
  -2 \\
  0
\end{pmatrix}\,.
\ee
The rotated 3D-index is given by
\begin{tiny}
\be
\begin{aligned}
\label{Irot237}
\Irot_{(-2,3,7)}(n,n')(q) =
&\sum_{k_1,k_2,k_3 \in \BZ} (-q^{1/2})^{k_1-2k_2-n+n^\prime}q^{k_3(n+n^\prime)/2}
\ID(k_1-2k_2-2k_3 +17n -17n^\prime,k_2+n-n^\prime)\\
&\hspace{-2cm} \times \ID(-k_1+k_2+n-n^\prime, k_1 -2k_2 - n + n^\prime)
\ID(2k_2 +n -n^\prime, k_1-2k_2-k_3+8n-8n^\prime)\,.
\end{aligned}
\ee
\end{tiny}

So, in our final experiment we computed the rotated 3D-index of the $(-2,3,7)$
pretzel-knot, and more precisely the $6 \times 6$ matrix $\Irot_{(-2,3,7)}(q)[6]$.
To give an idea of what this involves, the leading term of the above matrix is
\begin{tiny}
\be
\Irot_{(-2,3,7)}(q)[6] =
\begin{pmatrix}
1 & -q^{-9/2} & q^{-19} & -q^{-87/2} & q^{-78} & -q^{-245/2} \\
-q^{9/2} & 6 q^2 & -q^{-27/2} & q^{-38} & -q^{-145/2} & q^{-117} \\
q^{17} & -q^{27/2} & q & -q^{-45/2} & q^{57} & -q^{-203/2} \\
-q^{75/2} & q^{34} & -q^{45/2} & q^4 & -q^{-63/2} & q^{-76} \\
q^{66} & -q^{125/2} & q^{51} & -q^{63/2} & q^2 & -q^{-81/2} \\
-q^{205/2} & q^{99} & -q^{175/2} & q^{68} & -q^{81/2} & q^6
\end{pmatrix}
\ee
\end{tiny}

\noindent
and this alone required an internal truncation of the summand of~\eqref{Irot237}
up to $O(q^{103})$. For safety, we computed up to $O(q^{160})$ and we found out that
the last computed coefficients of $\Irot_{(-2,3,7)}(q)[6]$ were given by
{\fontsize{4.2}{4}\selectfont
\begin{equation*}
\begin{pmatrix}
3099301802486871 q^{158} & 15368338814987064 q^{315/2} & 
39577501827964202 q^{158} & -717771103116611523 q^{315/2} &
-7908419005020915850 q^{158} & 1907856058463675359575 q^{315/2}
\\
-2510483414752309 q^{315/2} & 3797180920247821 q^{158} & 
46280099948395184 q^{315/2} & -661349858819489021 q^{158} & 
6373738664932074312 q^{315/2} & 1164148757149541167314 q^{158}
\\
-830392595916755 q^{158} & -1589679235709546 q^{315/2} & 
5002197250330240 q^{158} & -59052244117713785 q^{315/2} & 
4279809698340893447 q^{158} & -25447538708964750026 q^{315/2}
\\
21883932028960 q^{315/2} & 52039830772006 q^{158} &
-208430252255007 q^{315/2} & 5021231467477637 q^{158} &
-203334247925102214 q^{315/2} & -14980307260595602909 q^{158}
\\
68212497673 q^{158} & -14703374329 q^{315/2} &
-986065940989 q^{158} & 1182082042782 q^{315/2} &
3294633659679268 q^{158} & 225454885754595400 q^{315/2}
\\
7690268 q^{315/2} & 27909767 q^{158} &
-486018210 q^{315/2} & -12829067397 q^{158} &
3046756706011 q^{315/2} & 1068804228132263 q^{158}
\end{pmatrix}
\end{equation*}
}
\noindent
On the other hand, the determinant of $\Irot_{(-2,3,7)}(q)[6]$ to that
precision was given by
\be
\det(\Irot_{(-2,3,7)}(q)[6])= q^{-15}(1 - q)^2 (1 - q^2)^4 (1 - q^3)^4 (1 - q^4)^2
+ O(q^{160}) \,.
\ee
But more reassuring was the fact that repeating the computation of
$\Irot_{(-2,3,7),\calO}(q)[6]$ for the insertion $\wh z_2$ (corresponding to the
second shape), we found out that the new matrix had equally big coefficients of
$q$-series, but the quotient
$$
Q_{(-2,3,7),\calO}(q)=\Irot_{(-2,3,7),\calO}(q)[6](\Irot_{(-2,3,7)}(q)[6])^{-1}
$$
had entries short rational functions
{\fontsize{6}{4}\selectfont
\label{237cancel}
\begin{multline*}
Q_{(-2,3,7),\calO}(q) + O(q^{160}) = \frac{1}{(1-q^3)(1-q^4)} \cdot
\\ 
\begin{pmatrix}
  \vspace{0.2cm}
0 & q^{-1/2}(q^2+1) &
q^{-19}(q^2-1)^2 (q^2+1) &
q^{-77/2}(-q^4-q^2-1) & q^{-74}(q^4-1)^2 & q^{-221/2}
\\ \vspace{0.2cm}
q^{15/2} (q^4+q^3+2 q^2+q+1) & (q-1)^2 (q^4+q^3+2 q^2+q+1) &
-q^{-23/2}(q+1) (q^2+1)^2 &
q^{-37}(q^2+1) (q^3-1)^2 & q^{-131/2}(q^2+1) & 0
\\ \vspace{0.2cm}
0 & q^{37/2} (q^2+1) & (q^2-1)^2 (q^2+1) &
q^{-39/2}(-q^4-q^2-1) & q^{-55}(q^4-1)^2 & q^{-183/2}
\\ \vspace{0.2cm}
q^{89/2} (q^4+q^3+2 q^2+q+1) &
(q-1)^2 q^{39} (q^4+q^3+2 q^2+q+1) &
-q^{51/2} (q+1) (q^2+1)^2 &
(q^2+1) (q^3-1)^2 & q^{-57/2}(q^2+1) & 0
\\ \vspace{0.2cm}
0 & q^{147/2} (q^2+1) &
q^{58} (q^2-1)^2 (q^2+1) & -q^{71/2} (q^4+q^2+1) &
(q^4-1)^2 & q^{-73/2}
\\ \vspace{0.2cm}
q^{235/2} (q^4+q^3+2 q^2+q+1) &
(q-1)^2 q^{112} (q^4+q^3+2 q^2+q+1) &
-q^{197/2} (q+1) (q^2+1)^2 &
q^{75} (q^2+1) (q^3-1)^2 & q^{89/2} (q^2+1) & 0
\end{pmatrix} 
\end{multline*}
}
\noindent
Surely this cancellation is not an accident, and it is a confirmation that our
computational method and Corollary~\ref{cor.Irotdesc} of Conjecture~\ref{conj.1}
are correct.

Incidentally, the $3 \times 3$ matrices $\Irot_{(-2,3,7),\calO}(q)[3]$
and $\Irot_{(-2,3,7),\calO}(q)[3]$ obey no rationality property similar to
Equation~\eqref{237cancel}, as one would not expect. 

\subsection*{Acknowledgements} 

The authors wish to thank Tudor Dimofte, Rinat Kashaev, Marcos Mari\~{n}o,
Campbell Wheeler and Don Zagier for many enlightening conversations. ZD would
like to thank International Center for Mathematics, SUSTech for hospitality
where this work was initiated. ZD is supported by KIAS individual Grant PG076902.


\appendix

\section{The holomorphic blocks of the $4_1$ knot}
\label{app.41}

The $4_1$ knot has two colored holomorphic blocks of the $4_1$ knot given by
$q$-hypergeometric formulas in~\cite[Prop.8]{GW:periods} as follows:
\be
\label{41h0}
h_{4_1,n}^{(0)}(q) = (-1)^nq^{|n|(2|n|+1)/2}
\sum_{k=0}^{\infty}(-1)^{k}\frac{q^{k(k+1)/2+|n|k}}{(q;q)_{k}(q;q)_{k+2|n|}}\,,
\ee
and
\be
\label{41h1}
\begin{small}
\begin{aligned}
h^{(1)}_{4_1,n}(q)
&=
(-1)^n q^{|n|(2|n|+1)/2}\sum_{k=0}^{\infty}\left(-4E_{1}(q)+
  \sum_{\ell=1}^{k+2|n|}\frac{1+q^\ell}{1-q^\ell}
  +\sum_{\ell=1}^{k}\frac{1+q^\ell}{1-q^\ell}\right)
(-1)^k\frac{q^{k(k+1)/2+|n|k}}{(q;q)_k(q,q)_{k+2|n|}}\\
&\quad-2(-1)^n q^{|n|(2|n|-1)/2}\sum_{k=0}^{2|n|-1}(-1)^k\frac{q^{k(k+1)/2-|n|k}
(q^{-1},q^{-1})_{2|n|-1-k}}{(q;q)_{k}}\,,
\end{aligned}
\end{small}
\ee
for $|q| \neq 1$. Here, for a positive integer $\ell$, we define
$
E_{\ell}(q)=\frac{\zeta(1-\ell)}{2}+\sum_{s=1}^{\infty}s^{\ell-1}\frac{q^{s}}{1-q^s}\,,
$
(where $\zeta(s)$ is the Riemann zeta function), analytic for $|q|<1$ and
extended to $|q|>1$ by the symmetry $E_{\ell}(q^{-1}) = -E_{\ell}(q)$.

\section{The holomorphic blocks of the $5_2$ knot}
\label{app.52}

The $5_2$ knot has three colored holomorphic blocks $h^{(\a)}_{5_2,n}(q)$ for $\a=0,1,2$.
They were given explicitly in~\cite[Lem.12]{GW:periods}, and we copy the
answer for the benefit of the reader. Using the $q$-harmonic functions
\be
H_{n}(q) = \sum_{j=1}^n \frac{q^j}{1-q^j}, \qquad
H^{(2)}_{n}(q) = \sum_{j=1}^n \frac{q^j}{(1-q^j)^2} 
\ee
we have:
\be
\label{52.h0}
\begin{aligned}
h_{5_2,n}^{(0)}(q)
&=
(-1)^nq^{|n|/2}\sum_{k=0}^{\infty}
\frac{q^{|n|k}}{(q^{-1};q^{-1})_{k}(q;q)_{k+2|n|}(q;q)_{k+|n|}}\,,
\end{aligned}
\ee
\be
\label{52.h1}
\begin{aligned}
h_{5_2,n}^{(1)}(q)
&=
-(-1)^nq^{|n|/2}\sum_{k=0}^{\infty}
\frac{q^{|n|k}}{(q;q)_{k+2|n|}(q^{-1};q^{-1})_{k}(q;q)_{k+|n|}}\\
&\qquad\times\left(k+|n|-\frac{1}{4}-3E_{1}(q)
  +H_{k}(q)
  +H_{k+|n|}(q)
  +H_{k+2|n|}(q)\right)\\
&\quad+q^{-n^2/2}\sum_{k=0}^{|n|-1}
\frac{(q^{-1},q^{-1})_{|n|-1-k}}{(q^{-1},q^{-1})_{k}(q;q)_{k+|n|}}\,,
\end{aligned}
\ee
and
\be
\label{52.h2}
\begin{aligned}
h_{5_2,n}^{(2)}(q)
&=
(-1)^nq^{|n|/2}\sum_{k=0}^{\infty}
\frac{q^{|n|k}}{(q^{-1};q^{-1})_{k}(q;q)_{k+|n|}(q;q)_{k+2|n|}}\\
&\qquad\times\Bigg(
E_{2}(q)+\frac{1}{8}
-H_{k}^{(2)}(q)
-H_{k+|n|}^{(2)}(q)
-H_{k+2|n|}^{(2)}(q)\\
&\qquad\qquad-\bigg(k+|n|-\frac{1}{4}-3E_{1}(q)
+H_{k}(q)
+H_{k+|n|}(q)
+H_{k+2|n|}(q)\bigg)^2\Bigg)\\
&+
2q^{-n^2/2}\sum_{k=0}^{|n|-1}
\frac{(q^{-1},q^{-1})_{|n|-1-k}}{(q^{-1},q^{-1})_{k}(q;q)_{k+|n|}}\\
&\qquad\times\Bigg(|n|-\frac{3}{4}-3E_{1}(q)
+H_{k}(q)
+H_{k+|n|}(q)
+H_{|n|-k-1}(q)\Bigg)\\
&-
2(-1)^{n}q^{-|n|/2}\sum_{k=0}^{|n|-1}q^{-|n|k}
\frac{(q^{-1};q^{-1})_{2|n|-k-1}(q^{-1};q^{-1})_{|n|-k-1}}{(q^{-1};q^{-1})_{k}}\,,
\end{aligned}
\ee
for $|q| \neq 1$.

\section{NZ matrices and the 3D-index}
\label{app.indexdef}

Since there are various formulas for the 3D-index in the literature, let us present
our conventions briefly.

Let $\calT$ be an ideal triangulation with $N$ tetrahedra of a 1-cusped hyperbolic
3-manifold $M$ equipped with a symplectic basis $\mu$ and $\lambda$ of
$H_1(\pt M,\BZ)$ and such that $\lambda$ is the homological longitude. Then
the edge gluing equations together with the peripheral equations are encoded by
three $(N+2)\times N$ matrices $\mb G$, $\mb G'$ and $\mb G''$ whose rows are indexed
by the edges, the meridian and the longitude and the columns indexed by tetrahedra.
The gluing equations in logarithmic form are given by
\be
\label{GGG}
\sum_{j=1}^N \big( \mb G_{ij} \log z_j + \mb G'_{ij} \log z'_j +
\mb G''_{ij} \log z''_j \big) = \pi i \,\, \boldsymbol\eta_i, \qquad i=1,\dots,N+2
\ee
where $\boldsymbol\eta=(2,\dots,2,0,0)^t \in \BZ^{N+2}$.

If we eliminate the variable $z'$ in each tetrahedron using $z z' z''=-1$,
we obtain the matrices $\mb A=\mb G-\mb G'$, $\mb B=\mb G''-\mb G'$ and the vector
$\boldsymbol \nu=(2,\dots,2,0,0)^t- \mb G'(1,\dots,1)^t$, and the gluing equations
take the form
\be
\label{ABnu}
\sum_{j=1}^N \big( \mb A_{ij} \log z_j + \mb B_{ij} \log z''_j \big)
= \pi i \,\, \boldsymbol\nu_i, \qquad i=1,\dots,N+2 \,.
\ee

Let $\mb a_j$ and $\mb b_j$ denote the $j$-th column of $\mb A$ and $\mb B$,
respectively. For integers $m$ and $e$, consider the vector
$\mb k = (k_1,\dots,k_{N-1},0,e,-m/2)$. Then, the 3D-index of~\cite{DGG1} is given
by~\cite{DGG1} (see also~\cite[Sec.4.5]{GHRS})


\be
\label{DGG}
I_\calT(m,e)(q) = \sum_{k_1,\dots,k_{N-1} \in \BZ}
(-q^{1/2})^{\boldsymbol \nu \cdot \mb k}
\prod_{j=1}^N \ID(-\mb b_j \cdot \mb k,\mb a_j \cdot \mb k)(q) 
\ee
and the rotated 3D-index is given by~\cite[Sec.2.1]{GW:periods}
\be
\label{Irotd1}
\Irot_\calT(n,n')(q) =\sum_{e \in \BZ} I_\calT(n-n',e)(q) q^{e(n+n')/2} \,.
\ee

Let us define the $N\times N$ matrices $A$ and $B$ obtained by removing the
$N$ and $N+2$ rows of $\mb A$ and $\mb B$, respectively. In other words, the
rows of $A$ and $B$ correspond to the first $N-1$ edge gluing equations and the
meridian gluing equation, respectively. Let $(\lambda_1,\dots,\lambda_N)$ and
$(\lambda''_1,\dots,\lambda''_N)$ denote \emph{half} the last row of
$\mb A$ and $\mb B$ respectively. We assume that these are vectors of integers and
this can be arranged by adding, if necessary, an integer multiple of some of the
first $N$ rows of $\mb A$ and $\mb B$ to the last row. Let $a_j$ and $b_j$ denote the
$j$-th column of $A$ and $B$, respectively, and let $k=(k_1,\dots,k_N)$.
Let $\nu \in \BZ^N$ be obtained from $\boldsymbol\nu \in \BZ^{N+2}$ by removing the
$N$-th and the $N+2$ entry of it, and let $\nu_\lambda$ denote half of the last entry
of $\boldsymbol\nu$. 

Then, combining~\eqref{DGG} and~\eqref{Irotd1} (where we rename its summation variable
from $e$ to $k_N$) we obtain that
\begin{small}
\be
\label{Irotd2}
\Irot_\calT(n,n')(q) =\sum_{k \in \BZ^{N}}
(-q^{1/2})^{\mb \nu \cdot k-(n-n')\nu_\lambda} q^{k_N(n+n')/2}
\prod_{j=1}^N \ID(\lambda''_j(n-n') -b_j \cdot k, -\lambda_j(n-n') + a_j \cdot k)(q)
\,.
\ee
\end{small}


\bibliographystyle{plain}
\bibliography{biblio}

\begin{thebibliography}{10}

\bibitem{Agarwal}
Prarit Agarwal, Dongmin Gang, Sangmin Lee, and Mauricio Romo.
\newblock Quantum trace map for 3-manifolds and a length conjecture.
\newblock Preprint 2022,
  \href{https://arxiv.org/abs/2203.15985}{arXiv:2203.15985}.

\bibitem{AK:TQFT}
J{\o}rgen~Ellegaard Andersen and Rinat Kashaev.
\newblock A {TQFT} from {Q}uantum {T}eichm\"uller theory.
\newblock {\em Comm. Math. Phys.}, 330(3):887--934, 2014.

\bibitem{Beem}
Christopher Beem, Tudor Dimofte, and Sara Pasquetti.
\newblock Holomorphic blocks in three dimensions.
\newblock {\em J. High Energy Phys.}, (12):177, front matter+118, 2014.

\bibitem{CCGLS}
Daryl Cooper, Marc Culler, Henry Gillet, Daryl Long, and Peter Shalen.
\newblock Plane curves associated to character varieties of {$3$}-manifolds.
\newblock {\em Invent. Math.}, 118(1):47--84, 1994.

\bibitem{Dimofte:riemann}
Tudor Dimofte.
\newblock Quantum {R}iemann surfaces in {C}hern-{S}imons theory.
\newblock {\em Adv. Theor. Math. Phys.}, 17(3):479--599, 2013.

\bibitem{Dimofte:3dsuper}
Tudor Dimofte.
\newblock 3d superconformal theories from three-manifolds.
\newblock In {\em New dualities of sypersymmetric gauge theories}, Math. Phys.
  Stud., pages 339--373. Springer, Cham, 2016.

\bibitem{DimofteGabellaGoncharov2016}
Tudor Dimofte, Maxime Gabella, and Alexander~B. Goncharov.
\newblock {K-Decompositions and 3d Gauge Theories}.
\newblock {\em JHEP}, 11:151, 2016.

\bibitem{DGG1}
Tudor Dimofte, Davide Gaiotto, and Sergei Gukov.
\newblock 3-manifolds and 3d indices.
\newblock {\em Adv. Theor. Math. Phys.}, 17(5):975--1076, 2013.

\bibitem{DGG2}
Tudor Dimofte, Davide Gaiotto, and Sergei Gukov.
\newblock Gauge theories labelled by three-manifolds.
\newblock {\em Comm. Math. Phys.}, 325(2):367--419, 2014.

\bibitem{DimofteGukovHollands2011}
Tudor Dimofte, Sergei Gukov, and Lotte Hollands.
\newblock {Vortex Counting and Lagrangian 3-manifolds}.
\newblock {\em Lett. Math. Phys.}, 98:225--287, 2011.

\bibitem{Gang:aspects}
Dongmin Gang, Nakwoo Kim, Mauricio Romo, and Masahito Yamazaki.
\newblock Aspects of defects in 3d-3d correspondence.
\newblock {\em J. High Energy Phys.}, (10):062, front matter+99, 2016.

\bibitem{Ga:AJ}
Stavros Garoufalidis.
\newblock On the characteristic and deformation varieties of a knot.
\newblock In {\em Proceedings of the {C}asson {F}est}, volume~7 of {\em Geom.
  Topol. Monogr.}, pages 291--309 (electronic). Geom. Topol. Publ., Coventry,
  2004.

\bibitem{GGM}
Stavros Garoufalidis, Jie Gu, and Marcos Mari\~{n}o.
\newblock The resurgent structure of quantum knot invariants.
\newblock {\em Comm. Math. Phys.}, 386(1):469--493, 2021.

\bibitem{GGMW:trivial}
Stavros Garoufalidis, Jie Gu, Marcos Mari\~no, and Campbell Wheeler.
\newblock Resurgence of {C}hern-{S}imons theory at the trivial flat connection.
\newblock Preprint 2021,
  \href{https://arxiv.org/abs/2111.04763}{arXiv:2111.04763}.

\bibitem{GHRS}
Stavros Garoufalidis, Craig Hodgson, Hyam Rubinstein, and Henry Segerman.
\newblock 1-efficient triangulations and the index of a cusped hyperbolic
  3-manifold.
\newblock {\em Geom. Topol.}, 19(5):2619--2689, 2015.

\bibitem{GK:desc}
Stavros Garoufalidis and Rinat Kashaev.
\newblock The descendant colored {J}ones polynomials.
\newblock Preprint 2021,
  \href{https://arxiv.org/abs/2108.07553}{arXiv:2108.07553}.

\bibitem{GK:qseries}
Stavros Garoufalidis and Rinat Kashaev.
\newblock From state integrals to {$q$}-series.
\newblock {\em Math. Res. Lett.}, 24(3):781--801, 2017.

\bibitem{GK:mero}
Stavros Garoufalidis and Rinat Kashaev.
\newblock A meromorphic extension of the 3{D} index.
\newblock {\em Res. Math. Sci.}, 6(1):Paper No. 8, 34, 2019.

\bibitem{GL:qholo}
Stavros Garoufalidis and Thang~T.Q. L{\^e}.
\newblock The colored {J}ones function is {$q$}-holonomic.
\newblock {\em Geom. Topol.}, 9:1253--1293 (electronic), 2005.

\bibitem{GW:periods}
Stavros Garoufalidis and Campbell Wheeler.
\newblock Periods, the meromorphic {3D}-index and the {T}uraev--{V}iro
  invariant.
\newblock Preprint 2022,
  \href{https://arxiv.org/abs/2209.02843}{arXiv:2209.02843}.

\bibitem{GZ:qseries}
Stavros Garoufalidis and Don Zagier.
\newblock Knots and their related $q$-series.
\newblock Preprint 2021.

\bibitem{GZ:kashaev}
Stavros Garoufalidis and Don Zagier.
\newblock Knots, perturbative series and quantum modularity.
\newblock Preprint 2021,
  \href{https://arxiv.org/abs/2111.06645}{arXiv:2111.06645}.

\bibitem{N:comb}
Walter Neumann.
\newblock Combinatorics of triangulations and the {C}hern-{S}imons invariant
  for hyperbolic {$3$}-manifolds.
\newblock In {\em Topology '90 ({C}olumbus, {OH}, 1990)}, volume~1 of {\em Ohio
  State Univ. Math. Res. Inst. Publ.}, pages 243--271. de Gruyter, Berlin,
  1992.

\bibitem{NZ}
Walter Neumann and Don Zagier.
\newblock Volumes of hyperbolic three-manifolds.
\newblock {\em Topology}, 24(3):307--332, 1985.

\bibitem{AB}
Marko Petkov\v{s}ek, Herbert~S. Wilf, and Doron Zeilberger.
\newblock {\em {$A=B$}}.
\newblock A K Peters, Ltd., Wellesley, MA, 1996.
\newblock With a foreword by Donald E. Knuth, With a separately available
  computer disk.

\bibitem{TerashimaYamazaki2011}
Yuji Terashima and Masahito Yamazaki.
\newblock {SL(2,R) Chern-Simons, Liouville, and Gauge Theory on Duality Walls}.
\newblock {\em JHEP}, 08:135, 2011.

\bibitem{WZ}
Herbert~S. Wilf and Doron Zeilberger.
\newblock An algorithmic proof theory for hypergeometric (ordinary and
  ``{$q$}'') multisum/integral identities.
\newblock {\em Invent. Math.}, 108(3):575--633, 1992.

\end{thebibliography}
\end{document}